# PERTURBATION OF ORTHOGONAL POLYNOMIALS ON AN ARC OF THE UNIT CIRCLE

Leonid Golinskii, Paul Nevai, and Walter Van Assche

Institute for Low Temperature Physics and Engineering,
Ohio State University, and Katholieke Universiteit Leuven

ABSTRACT. Orthogonal polynomials on the unit circle are completely determined by their reflection coefficients through the Szegő recurrences. We assume that the reflection coefficients converge to some complex number $a$ with $0 < |a| < 1$. The polynomials then live essentially on the arc $\{ e^{i\theta} : \alpha \leq \theta \leq 2\pi - \alpha \}$ where $\cos \frac{\alpha}{2} \stackrel{\text{def}}{=} \sqrt{1-|a|^2}$ with $\alpha \in (0,\pi)$. We analyze the orthogonal polynomials by comparing them with the orthogonal polynomials with constant reflection coefficients, which were studied earlier by Ya. L. Geronimus and N. I. Akhiezer. In particular, we show that under certain assumptions on the rate of convergence of the reflection coefficients the orthogonality measure will be absolutely continuous on the arc. In addition, we also prove the unit circle analogue of M. G. Krein's characterization of compactly supported nonnegative Borel measures on the real line whose support contains one single limit point in terms of the corresponding system of orthogonal polynomials.

## 1. Introduction

Orthogonal polynomials on the unit circle $\mathbb{T} \stackrel{\text{def}}{=} \{ z \in \mathbb{C} : |z| = 1 \}$ are defined by

$$\int_{\mathbb{T}} \varphi_n(\mu) \overline{\varphi_m(\mu)} \, d\mu = \delta_{m,n} , \qquad m,n = 0,1,2,\dots ,$$

1991 *Mathematics Subject Classification.* 42C05.

W.V.A. is a Senior Research Associate of the Belgian National Fund for Scientific Research. This material is based upon work supported by the International Science Foundation under Grant No. U9S000 (L.G.), by the National Science Foundation under Grants No. DMS–9024901 and No. DMS–940577 (P.N.), and by the Mathematics Research Institute of The Ohio State University (W.V.A.). The research started while W.V.A. was visiting The Ohio State University in July, 1993, and during L.G.'s and P.N.'s visit to the Technion, Haifa, in April, 1994.

Typeset by $\mathcal{A}_{\mathcal{M}}\mathcal{S}$-TEX

where
$$\varphi_n(\mu, z) = \kappa_n(\mu) z^n + \text{ lower degree terms}, \qquad \kappa_n(\mu) > 0,$$

and $\mu$ is a probability measure in $[0, 2\pi)$ with infinite support. Here and in what follows, for simplicity, we refer to $\mu$ as a probability measure on $\mathbb{T}$; if $f$ is a function on $\mathbb{T}$ then $\int_{\mathbb{T}} f \, d\mu \stackrel{\text{def}}{=} \int_0^{2\pi} f(z) \, d\mu(\theta)$ where $z = e^{i\theta}$, and, for instance, if $\mu'$ is the Radon-Nikodym derivative of $\mu$ then $\mu'(z) \stackrel{\text{def}}{=} \mu'(\theta)$ if $z = e^{i\theta}$ with $\theta \in [0, 2\pi)$. Standard references are books by Szegő [42], Freud [11], Grenander and Szegő [20], Geronimus [18], and two reasonably recent surveys [22] and [30].

These orthogonal polynomials satisfy the (Szegő) recurrences

$$\kappa_{n-1}(\mu)\varphi_n(\mu, z) = z\kappa_n(\mu)\varphi_{n-1}(\mu, z) + \varphi_n(\mu, 0)\varphi_{n-1}^*(\mu, z), \qquad n \in \mathbb{N}, \qquad (1)$$

$$\kappa_{n-1}(\mu)\varphi_n^*(\mu, z) = \kappa_n(\mu)\varphi_{n-1}^*(\mu, z) + z\overline{\varphi_n(\mu, 0)}\varphi_{n-1}(\mu, z), \qquad n \in \mathbb{N}, \qquad (2)$$

(cf. [42, formula (11.4.7), p. 293] and [18, formulas (1.2) and (1.2'), p. 6]) where the *reversed* $*$–polynomial of a polynomial $\rho_n$ of degree $n$ is given by $\rho_n^*(z) \stackrel{\text{def}}{=} z^n \overline{\rho_n}(1/z)$. The monic orthogonal polynomials are $\Phi_n(\mu) \stackrel{\text{def}}{=} \kappa_n^{-1}(\mu) \varphi_n(\mu)$ and then the Szegő recurrences take the form

$$\begin{pmatrix} \varphi_n(\mu, z) \\ \varphi_n^*(\mu, z) \end{pmatrix} = \frac{\kappa_n(\mu)}{\kappa_{n-1}(\mu)} \begin{pmatrix} z & \Phi_n(\mu, 0) \\ z\overline{\Phi_n(\mu, 0)} & 1 \end{pmatrix} \begin{pmatrix} \varphi_{n-1}(\mu, z) \\ \varphi_{n-1}^*(\mu, z) \end{pmatrix}, \qquad n \in \mathbb{N},$$

and

$$\begin{pmatrix} \Phi_n(\mu, z) \\ \Phi_n^*(\mu, z) \end{pmatrix} = \begin{pmatrix} z & \Phi_n(\mu, 0) \\ z\overline{\Phi_n(\mu, 0)} & 1 \end{pmatrix} \begin{pmatrix} \Phi_{n-1}(\mu, z) \\ \Phi_{n-1}^*(\mu, z) \end{pmatrix}, \qquad n \in \mathbb{N}.$$

The coefficients $\{\Phi_n(\mu, 0)\}_{n \in \mathbb{N}}$ describe completely not only the monic orthogonal polynomials on $\mathbb{T}$, but also the orthonormal polynomials since

$$\frac{\kappa_{n-1}^2(\mu)}{\kappa_n^2(\mu)} = 1 - |\Phi_n(\mu, 0)|^2 \quad \& \quad \kappa_n(\mu) = \left[\prod_{k=1}^n (1 - |\Phi_k(\mu, 0)|^2)\right]^{-\frac{1}{2}}, \qquad n \in \mathbb{N}, \qquad (3)$$

(cf. [18, p. 7]). The coefficients $\{\Phi_n(\mu, 0)\}_{n \in \mathbb{N}}$ are known as *reflection* and/or *recursion* and/or *Szegő* and/or *Schur* coefficients.[1]

Perhaps the easiest proof of (1) (and (2)) is using another set of (equivalent Szegő) recurrences which can be stated as

$$\kappa_n(\mu)\varphi_n^*(\mu, z) = \kappa_{n-1}(\mu)\varphi_{n-1}^*(\mu, z) + \overline{\varphi_n(\mu, 0)}\varphi_n(\mu, z), \qquad n \in \mathbb{N}, \qquad (4)$$

or, equivalently,

$$\kappa_n(\mu)\varphi_n^*(\mu, z) = \sum_{k=0}^n \overline{\varphi_k(\mu, 0)}\varphi_k(\mu, z), \qquad n \in \mathbb{N}. \qquad (5)$$

---

[1] Sometimes the sequence $\{-\overline{\Phi_{n+1}(\mu, 0)}\}_{n \in \mathbb{N}}$ is referred to by these names.

The latter follows directly from the fact that on the unit circle $\bar{z} = 1/z$, and, thus, two naturally arising extremal problems will have the same solutions (cf. [42, formula (11.4.7), p. 293] and [11, Theorem 5.1.8, p. 195]). If we replace $n$ by $n+1$ in (1) and then eliminate $\varphi_n^*(\mu)$ from it using (5), we obtain

$$z\varphi_n(\mu,z) = \frac{\kappa_n(\mu)}{\kappa_{n+1}(\mu)}\varphi_{n+1}(\mu,z) - \sum_{k=0}^{n} \frac{\kappa_k(\mu)}{\kappa_n(\mu)} \Phi_{n+1}(\mu,0)\, \overline{\Phi_k(\mu,0)}\, \varphi_k(\mu,z)\,, \qquad n \in \mathbb{N}\,, \quad (6)$$

which is the orthogonal Fourier expansion of $z\varphi_n(\mu,z)$ in terms of the orthogonal polynomial system $\{\varphi_n(\mu)\}_{n=0}^{\infty}$. The latter formula is useful when considering matrix representation of the multiplication operator in terms of $\{\varphi_n(\mu)\}_{n\in\mathbb{N}}$ (cf. Section 3).

Since all the zeros of $\Phi_n$ are inside the unit circle (cf. [18, Section 8, p. 9]), $|\Phi_n(\mu,0)| < 1$ for $n \in \mathbb{N}$. Conversely, when $\{\Phi_n(0)\}_{n\in\mathbb{N}}$ is a sequence of complex numbers with $|\Phi_n(0)| < 1$, then, by Favard's theorem on $\mathbb{T}$ (cf. [8] and [10]), the polynomials obtained by the Szegő recurrences are orthogonal with respect to a unique probability measure $\mu$ on $\mathbb{T}$ with infinite support so that $\Phi_n(0) = \Phi_n(\mu,0)$ for $n \in \mathbb{N}$.

Szegő's theory deals with the case when

$$\log \mu' \in L^1(\mathbb{T}) \iff \sum_{k=1}^{\infty} |\varphi_k(\mu,0)|^2 < \infty \iff \lim_{n\to\infty} \kappa_n(\mu) < \infty \iff$$

$$\overline{\mathbb{P}} \neq L^2(\mu,\mathbb{T}) \iff \prod_{k=1}^{\infty}(1 - |\Phi_k(\mu,0)|^2) > 0 \iff \sum_{k=1}^{\infty} |\Phi_k(\mu,0)|^2 < \infty \quad (7)$$

where $\mathbb{P}$ denotes the set of all algebraic polynomials with complex coefficients (cf. [42], [20], [18], or [15, Theorem VII, p. 751] and the references therein). Hence, under Szegő's condition, $\lim_{n\to\infty} \Phi_n(\mu,0) = 0$ always holds. The simplest case is when $\Phi_n(\mu,0) \equiv 0$ for $n \in \mathbb{N}$ which gives

$$\begin{pmatrix} \varphi_n(\mu,z) \\ \varphi_n^*(\mu,z) \end{pmatrix} = \begin{pmatrix} z & 0 \\ 0 & 1 \end{pmatrix} \begin{pmatrix} \varphi_{n-1}(\mu,z) \\ \varphi_{n-1}^*(\mu,z) \end{pmatrix} = \begin{pmatrix} z & 0 \\ 0 & 1 \end{pmatrix}^n \begin{pmatrix} 1 \\ 1 \end{pmatrix} = \begin{pmatrix} z^n \\ 1 \end{pmatrix}, \qquad n \in \mathbb{N}\,.$$

Given a probability measure $\mu$ on $\mathbb{T}$ with infinite support, in addition to $\{\varphi_n(\mu)\}_{n\in\mathbb{N}}$, we will also study polynomials of the *second kind* $\{\psi_n(\mu)\}_{n\in\mathbb{N}}$ which are defined as the orthogonal polynomials with reflection coefficients $\{-\Phi_n(\mu,0)\}_{n\in\mathbb{N}}$.[2] They can be computed simultaneously with the polynomials $\varphi_n(\mu)$ by Szegő type matrix recursions such as

$$\begin{pmatrix} \varphi_n(\mu,z) & \psi_n(\mu,z) \\ \varphi_n^*(\mu,z) & -\psi_n^*(\mu,z) \end{pmatrix} \qquad (8)$$

$$= \frac{\kappa_n(\mu)}{\kappa_{n-1}(\mu)} \begin{pmatrix} z & \Phi_n(\mu,0) \\ z\overline{\Phi_n(\mu,0)} & 1 \end{pmatrix} \begin{pmatrix} \varphi_{n-1}(\mu,z) & \psi_{n-1}(\mu,z) \\ \varphi_{n-1}^*(\mu,z) & -\psi_{n-1}^*(\mu,z) \end{pmatrix}.$$

---

[2]Note that the second kind orthogonal polynomials associated with the second kind orthogonal polynomials are the orthogonal polynomials one started with, as opposed to the case of orthogonal polynomials on the real line where this procedure leads to the notion of *associated polynomials*.

The advantage of using a matrix recurrence relation as opposed to a vector recursion is that one can manipulate inverse matrices whereas vectors are not invertible.

Since
$$\begin{vmatrix} \varphi_n(\mu,z) & \psi_n(\mu,z) \\ \varphi_n^*(\mu,z) & -\psi_n^*(\mu,z) \end{vmatrix} \equiv -2z^n, \qquad n \in \mathbb{N}, \tag{9}$$

(cf. [18, formula (1.17), p. 11]), we have
$$1 \leq |\varphi_n(\mu,z)| \, |\psi_n(\mu,z)|, \qquad z \in \mathbb{T}, \ n \in \mathbb{N}, \tag{10}$$

for every system of orthogonal polynomials on $\mathbb{T}$. In particular, if the orthogonal polynomials are bounded from above at a point on $\mathbb{T}$, then the second kind orthogonal polynomials are bounded from below at the same point, and vice versa. This leads us to the following simple but useful

**Lemma 1.** *Let $\mu$ be a probability measure on $\mathbb{T}$ having an infinite support. Assume that there is a closed circular arc $\Delta \subseteq \mathbb{T}$ and an infinite subset $\mathbb{N}^* \subseteq \mathbb{N}$ such that the second kind orthogonal polynomials $\{\psi_n(\mu)\}_{n \in \mathbb{N}^*}$ are uniformly bounded in $\Delta$. Then $\mu \in \mathrm{AC}(\Delta)$.*

*Proof of Lemma 1.* First, it follows from inequality (10) that the uniform boundedness of $\{\psi_n(\mu,z)\}_{n \in \mathbb{N}^*}$ in $\Delta$ implies the divergence of the series $\sum_{k \in \mathbb{N}} |\varphi_k(\mu,z)|^2$ in $\Delta$ so that $\mu$ has no mass points in $\Delta$ (cf. [15, Theorem VI, p. 750]). Second, for every such probability measure $\mu$,
$$\lim_{n \to \infty} \frac{1}{2\pi} \int_0^{2\pi} 1_{\Delta^*}(z) |\varphi_n(\mu,z)|^{-2} \, d\theta = \int_{\mathbb{T}} 1_{\Delta^*} \, d\mu, \qquad z = e^{i\theta}, \tag{11}$$

where $1_{\Delta^*}$ denotes the characteristic function of an arbitrary circular arc $\Delta^* \subseteq \Delta$ (cf. [27, Lemma 4.2, p. 248][3] ) so that by (10)
$$\mu(\Delta^*) \leq |\Delta^*| \sup_{n \in \mathbb{N}^*} \max_{z \in \Delta^*} |\psi_n(\mu,z)|^2, \qquad \Delta^* \subseteq \Delta,$$

where $|\Delta^*|$ denotes the normalized arc-length of $\Delta^*$, and, therefore, $\mu \in \mathrm{AC}(\Delta)$. $\square$

The following lemma is equally simple and equally useful.

**Lemma 2.** *Let $\mu$ be a probability measure on $\mathbb{T}$ having an infinite support. Then the inequality*
$$1/\mu'(\theta) \leq 2\pi \limsup_{n(\in \mathbb{N}^*) \to \infty} |\varphi_n(\mu,z)|^2, \qquad \text{for a.e. } z = e^{i\theta} \in \mathbb{T}, \tag{12}$$

*holds for each infinite $\mathbb{N}^* \subseteq \mathbb{N}$. In particular, if there is a closed circular arc $\Delta \subseteq \mathbb{T}$ and an infinite subset $\mathbb{N}^* \subseteq \mathbb{N}$ such that the orthogonal polynomials $\{\varphi_n(\mu)\}_{n \in \mathbb{N}^*}$ are uniformly bounded in $\Delta$, then $1/\mu' \in L^\infty(\Delta)$.*

*Proof of Lemma 2.* By [27, Lemma 4.2, p. 248] (cf. the footnote in the proof of Lemma 1) formula (11) holds if $1_{\Delta^*}$ is the characteristic function of an arbitrary circular arc $\Delta^*$ as long as the endpoints of $\Delta^*$ are not mass points of $\mu$. Hence, a simple application of Fatou's lemma completes the proof. □

---

[3]To be precise, in [27, Lemma 4.2, p. 248] the analogue of this is proved for continuous functions which then extends to all Riemann-Stieltjes integrable functions by one-sided approximation, and if $\mu$ has no mass points in $\Delta$ then $1_{\Delta^*}$ is Riemann-Stieltjes integrable for all $\Delta^* \subseteq \Delta$.

We are interested in orthogonal polynomials $\{\varphi_n(\mu)\}_{n\in\mathbb{N}}$ on the unit circle satisfying either $\lim_{n\to\infty} \Phi_n(\mu, 0) = a$ or the seemingly more general $\lim_{n\to\infty} \tau^n \Phi_n(\mu, 0) = a$ where $a \in \mathbb{C}$ with $0 < |a| < 1$ and $\tau \in \mathbb{T}$. We want to compare such orthogonal polynomials with the system of orthogonal polynomials with constant reflection coefficients $\{a\}_{n\in\mathbb{N}}$. The latter polynomials have been studied earlier by Geronimus and Akhiezer in, for instance, [19, p. 93], [14], [17, §4.3], and [1].

As far as $\lim_{n\to\infty} \Phi_n(\mu, 0) = 0$ goes, this not only follows from Szegő's condition $\log \mu' \in L^1(\mathbb{T})$ but it also holds whenever $\mu'$ is positive almost everywhere in $\mathbb{T}$. The latter was proved by Rakhmanov (cf. [37, Theorem 1, p. 105]). As a matter of fact,

$$\mu' > 0 \text{ a.e.} \iff \lim_{n\to\infty} \sup_{\ell \in \mathbb{N}} \int_0^{2\pi} \left| \frac{|\varphi_n(\mu, e^{it})|^2}{|\varphi_{n+\ell}(\mu, e^{it})|^2} - 1 \right| dt = 0$$

and

$$\lim_{n\to\infty} \Phi_n(\mu, 0) = 0 \iff \lim_{n\to\infty} \inf_{\ell \in \mathbb{N}} \int_0^{2\pi} \left| \frac{|\varphi_n(\mu, e^{it})|^2}{|\varphi_{n+\ell}(\mu, e^{it})|^2} - 1 \right| dt = 0$$

(cf. [26, Theorems 2 and 3, p. 64], [32, Theorem 1.1, p. 295], [33, Theorem 4, p. 325], and [23, Theorem B, p. 192]).

Finally, we mention

$$\sum_{n=1}^{\infty} |\Phi_n(\mu, 0)| < \infty \implies \mu \in \text{AC}(\mathbb{T}) \quad \& \quad (\mu')^{\pm 1} \in \text{C}(\mathbb{T})$$

(cf. [18, Theorem 8.5, p. 163]), and

$$\sum_{n=1}^{\infty} |\Phi_n(\mu, 0)| < \infty \iff$$

$$\mu \in \text{AC}(\mathbb{T}) \quad \& \quad (\mu')^{\pm 1} \in \text{C}(\mathbb{T}) \quad \& \quad \sum_{n=-\infty}^{\infty} \left| \int_0^{2\pi} e^{int} d\mu(t) \right| < \infty$$

(cf. [3, Theorem 8.1, p. 483]).

## 2. The case of constant reflection coefficients

Given $a \in \mathbb{C}$ with $0 < |a| < 1$, let $\mu_a$ denote the probability measure on $\mathbb{T}$ for which the corresponding orthogonal polynomials satisfy $\Phi_n(\mu_a, 0) = a$ for $n \in \mathbb{N}$. For sake of convenience, we will denote the corresponding orthogonal polynomials and second kind orthogonal polynomials by $\widehat{\varphi}_n$ and $\widehat{\psi}_n$, respectively, that is, $\widehat{\varphi}_n \stackrel{\text{def}}{=} \varphi_n(\mu_a)$ and $\widehat{\psi}_n \stackrel{\text{def}}{=} \psi_n(\mu_a)$ for $n \in \mathbb{N}$. In addition,

$$\cos \tfrac{\alpha}{2} \stackrel{\text{def}}{=} \sqrt{1 - |a|^2}, \qquad \alpha \in (0, \pi), \tag{13}$$

$$\Delta_\alpha \stackrel{\text{def}}{=} \{ e^{i\theta} : \alpha \leq \theta \leq 2\pi - \alpha \} \quad \text{and} \quad \Delta_\alpha^o \stackrel{\text{def}}{=} \{ e^{i\theta} : \alpha < \theta < 2\pi - \alpha \} \,. \quad (14)$$

With this notation, for $n \in \mathbb{N}$ the orthogonal polynomials with constant reflection coefficients $\{a\}_{n \in \mathbb{N}}$ are given by

$$\begin{pmatrix} \widehat{\varphi}_n(z) \\ \widehat{\varphi}_n^*(z) \end{pmatrix} = \frac{1}{\sqrt{1-|a|^2}} \begin{pmatrix} z & a \\ z\overline{a} & 1 \end{pmatrix} \begin{pmatrix} \widehat{\varphi}_{n-1}(z) \\ \widehat{\varphi}_{n-1}^*(z) \end{pmatrix} = (1-|a|^2)^{-n/2} \begin{pmatrix} z & a \\ z\overline{a} & 1 \end{pmatrix}^n \begin{pmatrix} 1 \\ 1 \end{pmatrix} \quad (15)$$

(cf. (3)). Similarly,

$$\begin{pmatrix} \widehat{\psi}_n(z) \\ \widehat{\psi}_n^*(z) \end{pmatrix} = \frac{1}{\sqrt{1-|a|^2}} \begin{pmatrix} z & -a \\ -z\overline{a} & 1 \end{pmatrix} \begin{pmatrix} \widehat{\psi}_{n-1}(z) \\ \widehat{\psi}_{n-1}^*(z) \end{pmatrix} = (1-|a|^2)^{-n/2} \begin{pmatrix} z & -a \\ -z\overline{a} & 1 \end{pmatrix}^n \begin{pmatrix} 1 \\ 1 \end{pmatrix} . \quad (16)$$

For $z \neq e^{\pm i\alpha}$, write

$$\begin{pmatrix} z & a \\ z\overline{a} & 1 \end{pmatrix} = V \begin{pmatrix} z_1 & 0 \\ 0 & z_2 \end{pmatrix} V^{-1}$$

where $z_1$ and $z_2$ are the eigenvalues of $\begin{pmatrix} z & a \\ z\overline{a} & 1 \end{pmatrix}$. Then

$$z_1 = \frac{z+1+\sqrt{(z-e^{i\alpha})(z-e^{-i\alpha})}}{2} \,, \qquad z_2 = \frac{z+1-\sqrt{(z-e^{i\alpha})(z-e^{-i\alpha})}}{2} \,,$$

where that branch of the square root is chosen for which

$$\lim_{z \to \infty} \frac{\sqrt{(z-e^{i\alpha})(z-e^{-i\alpha})}}{z} = 1 \,.$$

Then

$$\widehat{\varphi}_n(z) = \frac{Az_1^n + Bz_2^n}{(1-|a|^2)^{n/2}} \,, \qquad \widehat{\varphi}_n^*(z) = \frac{Cz_1^n + Dz_2^n}{(1-|a|^2)^{n/2}} \,, \qquad n = 0, 1, \ldots,$$

where $A$, $B$, $C$, and $D$ are functions of $z$ which do not depend on $n$. Setting $n=0$ and $n=1$ one can evaluate $A$, $B$, $C$, and $D$, which yields

$$\widehat{\varphi}_n(z) = \frac{1}{(1-|a|^2)^{n/2}} \left( (z+a) \frac{z_1^n - z_2^n}{z_1 - z_2} - z(1-|a|^2) \frac{z_1^{n-1} - z_2^{n-1}}{z_1 - z_2} \right) , \qquad n \in \mathbb{N}, \quad (17)$$

and

$$\widehat{\varphi}_n^*(z) = \frac{1}{(1-|a|^2)^{n/2}} \left( (1+\overline{a}z) \frac{z_1^n - z_2^n}{z_1 - z_2} - z(1-|a|^2) \frac{z_1^{n-1} - z_2^{n-1}}{z_1 - z_2} \right) , \qquad n \in \mathbb{N}. \quad (18)$$

Similarly,

$$\widehat{\psi}_n(z) = \frac{1}{(1-|a|^2)^{n/2}} \left( (z-a) \frac{z_1^n - z_2^n}{z_1 - z_2} - z(1-|a|^2) \frac{z_1^{n-1} - z_2^{n-1}}{z_1 - z_2} \right) , \qquad n \in \mathbb{N}, \quad (19)$$

and

$$\widehat{\psi}_n^*(z) = \frac{1}{(1-|a|^2)^{n/2}} \left( (1-\overline{a}z)\frac{z_1^n - z_2^n}{z_1 - z_2} - z(1-|a|^2)\frac{z_1^{n-1} - z_2^{n-1}}{z_1 - z_2} \right), \qquad n \in \mathbb{N}. \quad (20)$$

When $z = e^{i\theta}$ with $\alpha \leq \theta \leq 2\pi - \alpha$ then $|z_1| = |z_2| = \sqrt{1-|a|^2}$, and, hence, there is a constant $C > 0$ such that

$$|\widehat{\varphi}_n(z)| \leq C \min(n, v_\alpha(\theta)), \qquad z = e^{i\theta}, \ z \in \Delta_\alpha, \ n \in \mathbb{N}, \quad (21)$$

where

$$v_\alpha(\theta) \stackrel{\text{def}}{=} \sin\tfrac{\theta}{2} \, |\cos\alpha - \cos\theta|^{-\frac{1}{2}}, \quad (22)$$

and, therefore, for every $0 < \epsilon < \pi - \alpha$ there is a constant $C(\epsilon) > 0$ such that

$$|\widehat{\varphi}_n(z)| \leq C(\epsilon), \qquad z = e^{i\theta}, \ \alpha + \epsilon \leq \theta \leq 2\pi - \alpha - \epsilon, \ n \in \mathbb{N}. \quad (23)$$

In other words, given a fixed norm for $2 \times 2$ matrices, in view of (9), (21), and (23), there is a constant $C^* > 0$ such that

$$\left\| \begin{pmatrix} \widehat{\varphi}_n(z) & \widehat{\psi}_n(z) \\ \widehat{\varphi}_n^*(z) & -\widehat{\psi}_n^*(z) \end{pmatrix}^{\pm 1} \right\| \leq C^* \min(n, v_\alpha(\theta)),$$
$$z = e^{i\theta}, \ z \in \Delta_\alpha, \ n \in \mathbb{N}, \quad (24)$$

and for every $0 < \epsilon < \pi - \alpha$ there is a constant $C^*(\epsilon) > 0$ such that

$$\left\| \begin{pmatrix} \widehat{\varphi}_n(z) & \widehat{\psi}_n(z) \\ \widehat{\varphi}_n^*(z) & -\widehat{\psi}_n^*(z) \end{pmatrix}^{\pm 1} \right\| \leq C^*(\epsilon), \qquad z = e^{i\theta}, \ \alpha + \epsilon \leq \theta \leq 2\pi - \alpha - \epsilon, \ n \in \mathbb{N}. \quad (25)$$

For these polynomials we have

$$\int_\alpha^{2\pi-\alpha} \widehat{\varphi}_n(z)\overline{\widehat{\varphi}_m(z)} \, \frac{\sqrt{\sin\frac{\theta+\alpha}{2}\sin\frac{\theta-\alpha}{2}}}{2\pi \sin\frac{\theta-\beta}{2}} \, d\theta + j_\beta \widehat{\varphi}_n(e^{i\beta})\overline{\widehat{\varphi}_m(e^{i\beta})} = \delta_{m,n}, \quad z = e^{i\theta}, \ m,n \in \mathbb{Z}^+,$$

with

$$e^{i\beta} = \frac{1-a}{1-\overline{a}}, \qquad j_\beta = \begin{cases} \frac{2|a|^2 - a - \overline{a}}{|1-a|}, & \text{if } |1-2a| > 1, \\ 0, & \text{if } |1-2a| \leq 1, \end{cases}$$

(cf. [19, formulas (XI.26) and (XI.27), p. 94]).[4]

---

[4] N.B. that our $a$ corresponds to $-\overline{a}$ in [19] (see footnote 2).

## 3. Spectral analysis

For orthogonal polynomials on the real line there is an intimate relationship with infinite Jacobi matrices containing the coefficients of the three-term recurrence relation for the orthonormal polynomials. These Jacobi matrices are symmetric tridiagonal matrices which can be extended to self-adjoint operators acting on the Hilbert space $\ell_2(\mathbb{N})$. Applying results from perturbation theory of self-adjoint operators then allow the interpretation of spectral properties of the Jacobi matrix as properties of the orthogonality measure. Such properties include (but are not limited to) information about its support and enables one to study its absolute continuity.

For orthogonal polynomials on the unit circle there is a similar relationship with infinite matrices, but instead of a self-adjoint tridiagonal matrix (for determinate moment problems on the real line) one deals with a unitary Hessenberg matrix (for measures outside the Szegő class).

Given a probability measure $\mu$ on $\mathbb{T}$ with infinite support, let $L^2(\mu, \mathbb{T})$ denote the Hilbert space of measurable, square-integrable functions on the unit circle $\mathbb{T}$ with the scalar product and norm

$$\langle f, g \rangle_\mu \stackrel{\text{def}}{=} \int_\mathbb{T} f \overline{g} \, d\mu \qquad \text{and} \qquad \|f\|_\mu \stackrel{\text{def}}{=} \sqrt{\langle f, f \rangle_\mu}.$$

By a theorem of A. N. Kolmogorov and M. G. Krein, the system $\{\varphi_n(\mu)\}_{n=0}^\infty$ forms an orthonormal basis in $L^2(\mu, \mathbb{T})$ if and only if $\log \mu' \notin L^1(\mathbb{T})$ (cf. (7) and [20, Theorem 3.3(a), p. 49]). All measures $\mu$ considered in this section are such that $\log \mu' \notin L^1(\mathbb{T})$.

In our investigations a key role is played by the unitary multiplication operator $U(\mu) : L^2(\mu, \mathbb{T}) \to L^2(\mu, \mathbb{T})$ defined by

$$[U(\mu) f](t) = t\, f(t), \qquad t \in \mathbb{T},\ f \in L^2(\mu, \mathbb{T}), \tag{26}$$

and its matrix representation $\widehat{U}(\mu)$ in the orthonormal basis $\{\varphi_n(\mu)\}_{n=0}^\infty$.

By (3) and (6) we have

$$\widehat{U}(\mu) = \begin{pmatrix} u_{00} & u_{01} & \cdots \\ u_{10} & u_{11} & \cdots \\ \vdots & \vdots & \ddots \end{pmatrix}, \qquad u_{ij} = \langle U(\mu)\varphi_j(\mu), \varphi_i(\mu) \rangle_\mu, \tag{27}$$

where

$$u_{ij} = \begin{cases} -\Phi_{j+1}(\mu, 0)\overline{\Phi_i(\mu, 0)} \prod_{k=i+1}^{j} \left(1 - |\Phi_k(\mu, 0)|^2\right)^{1/2}, & i < j+1, \\ \left(1 - |\Phi_{j+1}(\mu, 0)|^2\right)^{1/2}, & i = j+1, \\ 0, & i > j+1. \end{cases} \tag{28}$$

for $i, j = 0, 1, \ldots$ (cf. [43, formula (3), p. 409] and for a doubly infinite analogue see [13]). Infinite matrices such as (27)–(28) in which all entries below the subdiagonal vanish are called (upper) *Hessenberg* matrices.

We can view the infinite matrix (27)–(28) as a unitary operator $\widehat{U}(\mu) : \ell_2(\mathbb{N}) \to \ell_2(\mathbb{N})$ which is unitarily equivalent to the multiplication operator $U(\mu)$. In particular, $\operatorname{supp}(\mu)$ is equal to the spectrum of $\widehat{U}(\mu)$.

**Theorem 3.** (Geronimus, 1941) *Let $\mu$ be a probability measure on $\mathbb{T}$ with infinite support, and let $\{\varphi_n(\mu)\}_{n=0}^{\infty}$ be the corresponding orthogonal polynomials. Suppose $\lim_{n\to\infty} \Phi_n(\mu,0) = a$, where $a \in \mathbb{C}$ with $0 < |a| < 1$. Let $\cos\frac{\alpha}{2} \stackrel{\text{def}}{=} \sqrt{1-|a|^2}$ where $\alpha \in (0,\pi)$ and $\Delta_\gamma \stackrel{\text{def}}{=} \{e^{i\theta} : \gamma \leq \theta \leq 2\pi - \gamma\}$. Then $\Delta_\alpha \subseteq \mathrm{supp}(\mu)$ and $\mathrm{supp}(\mu) \setminus \Delta_\beta$ is finite for every $0 < \beta < \alpha$.*

*Remark 4.* Geronimus [14, Theorem 1', p. 205] used continued fractions to prove Theorem 3, whereas our proof is based on the spectral theory of unitary operators.[5]

*Remark 5.* The statement that "$\mathrm{supp}(\mu) \setminus \Delta_\beta$ is finite for every $0 < \beta < \alpha$" is just another way of saying that $\mathrm{supp}(\mu) \setminus \Delta_\alpha$ is at most countable whose limit points (if any) must belong to $\Delta_\alpha$.

*Proof of Theorem 3.* Along with the multiplication operator $U(\mu)$ (cf. (26)) consider the multiplication operator $U(\mu_a)$ on the space $L^2(\mu_a, \mathbb{T})$ where the measure $\mu_a$ corresponds to the constant reflection coefficients $\Phi_n(\mu_a, 0) = a$ for $n \in \mathbb{N}$ (see Section 2). To study both these operators simultaneously consider their matrix representations $\widehat{U}(\mu)$ (cf. (27)–(28)) and

$$\widehat{U}(\mu_a) = \begin{pmatrix} u_{00} & u_{01} & \cdots \\ u_{10} & u_{11} & \cdots \\ \vdots & \vdots & \ddots \end{pmatrix}, \qquad u_{ij} = \langle U(\mu_a)\varphi_j(\mu_a), \varphi_i(\mu_a)\rangle_{\mu_a},$$

where

$$u_{ij} = \begin{cases} -a\left(1-|a|^2\right)^{j/2}, & i = 0, \\ -|a|^2\left(1-|a|^2\right)^{(j-i)/2}, & 0 < i < j+1, \\ \left(1-|a|^2\right)^{1/2}, & i = j+1, \\ 0, & i > j+1. \end{cases}$$

for $i,j = 0,1,\ldots$. This follows from (27)–(28) applied with $\mu$ replaced by $\mu_a$. Both $\widehat{U}(\mu)$ and $\widehat{U}(\mu_a)$ are unitary operators acting on the *same* Hilbert space $\ell_2(\mathbb{N})$. As discussed in Section 2, the continuous spectrum of $\widehat{U}(\mu_a)$ is the arc $\Delta_\alpha$, and, in addition, $\widehat{U}(\mu_a)$ may have at most one eigenvalue which must be in $\mathbb{T}$ but is located outside this arc. Let $S : \ell_2(\mathbb{N}) \to \ell_2(\mathbb{N})$ be the shift operator given by the matrix representation

$$S = \begin{pmatrix} 0 & 1 & 0 & 0 & \cdots \\ 0 & 0 & 1 & 0 & \cdots \\ 0 & 0 & 0 & 1 & \cdots \\ 0 & 0 & 0 & 0 & \cdots \\ \vdots & \vdots & \vdots & \ddots & \ddots \end{pmatrix}.$$

Then we can write $\widehat{U}(\mu)$ as

$$\widehat{U}(\mu) = S^* D_{-1}(\mu) + \sum_{j=0}^{\infty} D_j(\mu) S^j \qquad (29)$$

---

[5]N.B. Instead of $\Delta_\alpha \subseteq \mathrm{supp}(\mu)$, Geronimus claims $\Delta_\alpha \subseteq \overline{\mathrm{supp}(\mu)}$, but, of course, $\mathrm{supp}(\mu)$ is closed.

where $S^*$ is the adjoint of $S$ and each $D_j(\mu) : \ell_2(\mathbb{N}) \to \ell_2(\mathbb{N})$ is a diagonal matrix whose main diagonal is equal to the $j$th diagonal above the main diagonal of $\widehat{U}(\mu)$, that is,

$$\operatorname{diag} D_j(\mu) = \begin{cases} \left\{ -\Phi_{j+i+1}(\mu,0)\overline{\Phi_i(\mu,0)} \prod_{k=i+1}^{i+j} \sqrt{1-|\Phi_k(\mu,0)|^2} \right\}_{i=0}^{\infty}, & j \in \mathbb{Z}^+ \\ \left\{ \sqrt{1-|\Phi_i(\mu,0)|^2} \right\}_{i=0}^{\infty}, & j = -1, \end{cases}$$

This infinite series representation for $\widehat{U}(\mu)$ converges in the operator norm. To see this note that $\|S\| = 1$ and

$$\begin{aligned} \|D_j(\mu)\| &= \sup_{i \geq 0} \left| \Phi_{j+i+1}(\mu,0)\overline{\Phi_i(\mu,0)} \right| \prod_{k=i+1}^{i+j} \sqrt{1-|\Phi_k(\mu,0)|^2} \\ &\leq \sup_{i \geq 0} \prod_{k=i+1}^{i+j} \sqrt{1-|\Phi_k(\mu,0)|^2}, \end{aligned} \tag{30}$$

for $j \in \mathbb{Z}^+$ since $|\Phi_n(\mu,0)| < 1$ for $n \in \mathbb{N}$. If $\lim_{n\to\infty} \Phi_n(\mu,0) = a$ with $a \neq 0$ then $\{\|D_j(\mu)\|\}_{j=-1}^{\infty}$ decreases exponentially, and, therefore the series in (29) converges uniformly.

Similarly, we can write $\widehat{U}(\mu_a)$ as a uniformly convergent series

$$\widehat{U}(\mu_a) = S^* D_{-1}(\mu_a) + \sum_{j=0}^{\infty} D_j(\mu_a) S^j.$$

Now consider the difference $\widehat{U}(\mu) - \widehat{U}(\mu_a)$. Given $j = -1, 0, 1, \ldots$, the difference $D_j(\mu) - D_j(\mu_a)$ is a compact operator since it is a diagonal operator for which the entries converge to zero. Linear combinations of compact operators and the product of a compact operator with a bounded operator remain compact. Moreover, the set of compact operators on a Hilbert space is closed. Hence, $\widehat{U}(\mu) - \widehat{U}(\mu_a)$ is a compact operator. But then H. Weyl's theorem[6] says that, except for eigenvalues, the spectra of $\widehat{U}(\mu)$ and $\widehat{U}(\mu_a)$ are the same. Since the spectra of $\widehat{U}(\mu)$ and $\widehat{U}(\mu_a)$ are equal to the supports of $\mu$ and $\mu_a$, respectively, this is precisely what had to be proved in view of the constant reflection coefficient case discussed in Section 2. □

We finish this section by proving the unit circle analogue of M. G. Krein's characterization of compactly supported nonnegative Borel measures on the real line whose support contains one single limit point, in terms of the corresponding system of orthogonal polynomials (cf. [2, Theorem 3, p. 231] and [4, Theorem 4.3.5, p. 117]). The general case which deals with finitely many limit points just as M. G. Krein dealt with measures on the real line (cf. [2, Theorem 2, p. 230] and [4, Theorem 4.6.2, p. 137]) remains open.

---

[6]H. Weyl proved his theorem for self-adjoint operators. We are using a more general version as given in [21, Problem 143, p. 91] and [44, Proposition 1, p. 62].

**Theorem 6.** *Let $\mu$ be a probability measure on $\mathbb{T}$ having an infinite support, and let $\tau \in \mathbb{T}$. Then the following statements are equivalent.*

*(i) The derived set of the support of $\mu$ is equal to $\{\tau\}$.*

*(ii) We have $\lim_{n\to\infty} \int_{\mathbb{T}} f |\varphi_n(\mu)|^2 \, d\mu = f(\tau)$ for every $\mu$-measurable function $f$ which is bounded on $\mathbb{T}$ and continuous at $\tau$.*

*(iii) We have $\lim_{n\to\infty} \int_{\mathbb{T}} f \varphi_n(\mu) \overline{\varphi_{n+k}(\mu)} \, d\mu = f(\tau) \delta_{0,k}$ for all $k \in \mathbb{Z}$ and uniformly in $k \in \mathbb{Z}^+$ for every $\mu$-measurable function $f$ which is bounded on $\mathbb{T}$ and continuous at $\tau$.*

*(iv) We have $\lim_{n\to\infty} \int_{\mathbb{T}} z |\varphi_n(\mu, z)|^2 \, d\mu(\theta) = \tau$ where $z = e^{i\theta}$.*

*(v) We have $\lim_{n\to\infty} \int_{\mathbb{T}} z \varphi_n(\mu, z) \overline{\varphi_{n+k}(\mu, z)} \, d\mu(\theta) = \tau \delta_{0,k}$ for all $k \in \mathbb{Z}$ and uniformly in $k \in \mathbb{Z}^+$ where $z = e^{i\theta}$.*

*(vi) We have $\lim_{n\to\infty} \Phi_{n+1}(\mu, 0) \overline{\Phi_n(\mu, 0)} = -\tau$.*

*Proof of Theorem 6.* For sake of brevity, $e^{i\theta}$ is denoted by $z$ in all integrals.

(i) $\Rightarrow$ (ii): Fix $\epsilon > 0$. Then, since $\varphi_n$'s are orthonormal,

$$\left| \int_{\mathbb{T}} f \, |\varphi_n(\mu)|^2 \, d\mu - f(\tau) \right|$$
$$= \left| \int_{\mathbb{T}} [f(z) - f(\tau)] \, |\varphi_n(\mu, z)|^2 \, d\mu(\theta) \right| \leq \int_{\mathbb{T}} |f(z) - f(\tau)| \, |\varphi_n(\mu, z)|^2 \, d\mu(\theta)$$
$$= \int_{|z-\tau|<\epsilon} |f(z) - f(\tau)| \, |\varphi_n(\mu, z)|^2 \, d\mu(\theta) + \int_{|z-\tau|\geq\epsilon} |f(z) - f(\tau)| \, |\varphi_n(\mu, z)|^2 \, d\mu(\theta)$$
$$\leq \sup_{|z-\tau|<\epsilon} |f(z) - f(\tau)| + \int_{|z-\tau|\geq\epsilon} |f(z) - f(\tau)| \, |\varphi_n(\mu, z)|^2 \, d\mu(\theta) \, .$$

It is well known that for every $z = e^{i\theta} \in \mathbb{T}$ we have $\sum_{n=0}^{\infty} |\varphi_n(\mu, z)|^2 \leq [\mu(\{\theta\})]^{-1}$. This follows, for instance, from the extremal property satisfied by the Christoffel functions (cf. [42, Theorem 11.3.1, p. 290]).[7] By (i), the set $\Delta_\epsilon \stackrel{\text{def}}{=} \text{supp}(\mu) \cap \{z \in \mathbb{T} : |z - \tau| \geq \epsilon\}$ is a finite collection of mass points of $\mu$ so that $\sum_{n=0}^{\infty} |\varphi_n(\mu, z)|^2$ converges uniformly on $\Delta_\epsilon$. In particular, $\lim_{n\to\infty} |\varphi_n(z)| = 0$ uniformly on $\Delta_\epsilon$ so that $\lim_{n\to\infty} \int_{|z-\tau|\geq\epsilon} |f(z) - f(\tau)| |\varphi_n(\mu, z)|^2 \, d\mu = 0$. Therefore,

$$\limsup_{n\to\infty} \left| \int_{\mathbb{T}} f |\varphi_n(\mu)|^2 \, d\mu - f(\tau) \right| \leq \sup_{|z-\tau|<\epsilon} |f(z) - f(\tau)| \, .$$

Letting $\epsilon \to 0$ yields (ii).

(ii) $\iff$ (iii): (ii) is a special case of (iii) whereas "(ii) $\Rightarrow$ (iii)" is proved as follows. Given $k \in \mathbb{Z}$ and $n \geq -k$, by orthogonality we have

$$\int_{\mathbb{T}} f \varphi_n(\mu) \overline{\varphi_{n+k}(\mu)} \, d\mu - f(\tau) \delta_{0,k} = \int_{\mathbb{T}} [f(z) - f(\tau)] \varphi_n(\mu, z) \overline{\varphi_{n+k}(\mu, z)} \, d\mu(\theta)$$

---

[7]Actually, we have $\sum_{n=0}^{\infty} |\varphi_n(\mu, z)|^2 = [\mu(\{\theta\})]^{-1}$, but the latter is a bit harder to prove as opposed to the inequality which is straightforward (cf. formula (7) on p. 453 and its proof on p. 444–445 in *Szegő's extremum problem on the unit circle* by A. Máté P. Nevai and V. Totik in Ann. Math. **134** (1991), 433–453).

so that Cauchy-Schwarz's inequality yields

$$\left| \int_{\mathbb{T}} f \varphi_n(\mu) \overline{\varphi_{n+k}(\mu)} \, d\mu - f(\tau) \delta_{0,k} \right|^2 \leq \int_{\mathbb{T}} |f(z) - f(\tau)|^2 \, |\varphi_n(\mu, z)|^2 \, d\mu(\theta) \, .$$

Therefore, (iii) follows directly from (ii).

(iii) $\Rightarrow$ (v): This is straightforward by letting $f(z) \stackrel{\text{def}}{=} z$.

(iv) $\iff$ (v): (iv) is a special case of (v) whereas "(iv) $\Rightarrow$ (v)" is proved as follows. Given $k \in \mathbb{Z}$ and $n \geq -k$, by orthogonality we have

$$\int_{\mathbb{T}} z \varphi_n(\mu, z) \overline{\varphi_{n+k}(\mu, z)} \, d\mu(\theta) - \tau \delta_{0,k} = \int_{\mathbb{T}} (z - \tau) \varphi_n(\mu, z) \overline{\varphi_{n+k}(\mu, z)} \, d\mu(\theta) \, .$$

Hence, Cauchy-Schwarz's inequality combined with (iv) and $\varphi_n$'s orthonormality yields

$$\left| \int_{\mathbb{T}} \varphi_n(\mu, z) \overline{\varphi_{n+k}(\mu, z)} \, d\mu(\theta) - \tau \delta_{0,k} \right|^2 \leq$$

$$\int_{\mathbb{T}} |z - \tau|^2 |\varphi_n(\mu, z)|^2 \, d\mu(\theta) = 2 - 2\Re\left\{ \overline{\tau} \int_{\mathbb{T}} z |\varphi_n(\mu, z)|^2 \, d\mu(\theta) \right\} \xrightarrow[n \to \infty]{} 0$$

so that (v) holds.

(iv) $\iff$ (vi): This follows by observing that, according to (6), $\int_{\mathbb{T}} z |\varphi_n(\mu, z)|^2 \, d\mu = -\Phi_{n+1}(\mu, 0) \overline{\Phi_n(\mu, 0)}$.

(v) $\Rightarrow$ (i): Consider $U(\mu)$ defined by (26) and it matrix representation $\widehat{U}(\mu)$ given by (27)–(28). Just like in the proof of Theorem 3 and borrowing the notation from there, write $\widehat{U}(\mu)$ as

$$\widehat{U}(\mu) = S^* D_{-1}(\mu) + \sum_{j=0}^{\infty} D_j(\mu) S^j \tag{31}$$

where $D_j(\mu)$ is a diagonal matrix whose diagonal is equal to the $j$th diagonal above the main diagonal of $\widehat{U}(\mu)$, that is,

$$\operatorname{diag} D_j(\mu) = \langle U(\mu) \varphi_{j+n}(\mu), \varphi_n(\mu) \rangle_{\mu} = \left\{ \int_0^{2\pi} z \varphi_{j+n}(\mu, z) \overline{\varphi_n(\mu, z)} \, d\mu(\theta) \right\}_{n=0}^{\infty}$$

for $j = -1, 0, 1, \ldots$.

Thus, the entries on the main diagonal of $D_0(\mu) - \tau \mathbb{I}$ and, for each $j \neq 0$, the entries on the main diagonal of $D_j(\mu)$ converge to 0 whenever (v) holds, that is, $D_0(\mu) - \tau \mathbb{I}$ and each such $D_j(\mu)$ are compact operators.[8]

By (vi) (which is equivalent to (v)) and

$$|\Phi_{n+1}(\mu, 0) \overline{\Phi_n(\mu, 0)}| \leq |\Phi_n(\mu, 0)| < 1 \, , \qquad n \in \mathbb{N} \, ,$$

---

[8] $\mathbb{I}$ denotes the identity matrix.

we have $\lim_{n\to\infty} |\Phi_n(\mu, 0)| = 1$. Thus, by (30), $\{\|D_j(\mu)\|\}_{j=-1}^{\infty}$ decreases exponentially, and then the series in (31) converges in the operator norm. Hence, $\widehat{U}(\mu) - \tau\mathbb{I} : \ell_2(\mathbb{N}) \to \ell_2(\mathbb{N})$ is a compact operator.

According to a theorem of F. Riesz and J. Schauder (cf. [40, Section 79, p. 187] and [39, Theorem VI.15, p. 203]) the spectrum of the operator $\widehat{U}(\mu) - \tau\mathbb{I}$ is an at most countable set with no limit point except possibly zero. Since $\mathrm{supp}(\mu)$ is an infinite set, the spectrum of the unitary operator $\widehat{U}(\mu) = \left(\widehat{U}(\mu) - \tau\mathbb{I}\right) + \tau\mathbb{I}$ is a countable set, it lies on the unit circle, and its only limit point is $\{\tau\}$. This is precisely what is stated in (i) since the spectrum of $\widehat{U}(\mu)$ is equal to $\mathrm{supp}(\mu)$. □

*Remark 7.* Geronimus points out that the one-point version of M. G. Krein's theorem (which was known to Stieltjes [41, note on p. 564–570]) may be paraphrased as follows (cf. [16, Theorem 32.2, p. 72]). If the sequence $\{\Phi_n(\mu, 0)\}_{n\in\mathbb{N}}$ is real then the derived set of the support of $\mu$ is equal to $\{-1\}$ whenever

$$\lim_{n\to\infty} (1 + \Phi_n(\mu, 0))(1 - \Phi_{n+1}(\mu, 0)) = 0.$$

It is easy to see that for real sequences $\{\Phi_n(\mu, 0)\}_{n\in\mathbb{N}}$ Geronimus' condition holds if and only if $\lim_{n\to\infty} \Phi_{n+1}(\mu, 0)\Phi_n(\mu, 0) = 1$.

As a matter of fact, for real sequences $\{a_n \in (-1, 1)\}_{n\in\mathbb{N}}$ or even for $\{a_n \in [-1, 1]\}_{n\in\mathbb{N}}$, the conditions (i) $\lim_{n\to\infty} a_n = a$ where $a$ is either $1$ or $-1$, (ii) $\lim_{n\to\infty} a_n a_{n+1} = 1$, and (iii) $\lim_{n\to\infty}(1 + a_n)(1 - a_{n+1}) = 0$, are all equivalent. The proof of this is an elementary exercise about $\liminf$'s and $\limsup$'s. For instance, if (iii) holds and $\limsup_{n\to\infty} a_n = 1$ then there is an infinite subsequence of the $a_n$'s which are all positive and then (iii) via induction implies that the $a_n$'s are all positive for sufficiently large $n$'s so that by (iii) we have $\lim_{n\to\infty} a_n = 1$. On the other hand, if (iii) holds and $\limsup_{n\to\infty} a_n < 1$, then $(1 + a_n)$ can be estimated from above by a constant multiple of $(1 + a_n)(1 - a_{n+1})$ for sufficiently large $n$'s so that by (iii) we have $\lim_{n\to\infty} a_n = -1$. Proving (i) from (ii) is even simpler.

*Remark 8.* If $\limsup_{n\to\infty} |\Phi_n(\mu, 0)| = 1$ then $\mu$ is *singular* (cf. [37, Lemma 4, p. 110] and for its operator theoretic analogue see [6]).[9]

*Remark 9.* Observe that $\Phi_n(\mu, 0) \stackrel{\mathrm{def}}{=} (1 - n^{-1})\exp(i\log n)$, $n \in \mathbb{N}$, is a divergent sequence for which (vi) in Theorem 6 holds.

**Example 10.** A. Zhedanov's example where, given $0 < q < 1$, $\Phi_n(\mu, 0) \stackrel{\mathrm{def}}{=} 2q^n - 1$, $n \in \mathbb{N}$, yields yet another sequence for which (vi) in Theorem 6 holds.

The details of A. Zhedanov's example are as follows. We construct a measure $\mu$ on $\mathbb{T}$ with $-1$ as the only limit point in its support and calculate its recursion coefficients explicitly.[10]

---

[9] We point out a typographical error in [37] appearing both in the Russian original and in the English translation. Namely, the condition $\limsup_{n\to\infty} |\Phi_n(\mu, 0)| = 1$ is stated as $\limsup_{n\to\infty} \Phi_n(\mu, 0) = 1$.

[10] We thank A. Zhedanov for his permission to include his example here (cf. [46, p. 1–2]).

Given $0 < q < 1$, consider the so-called (monic) "discrete $q$-Hermite polynomials" (also known as the Al-Salam & Carlitz polynomials)

$$H_n(x;q) = \sum_{k=0}^{[n/2]} \frac{(q;q)_n}{(q^2;q^2)_k (q;q)_{n-2k}} (-1)^k q^{k(k-1)} x^{n-2k}, \qquad 0 < q < 1,$$

where $(a;q)_n \stackrel{\text{def}}{=} \prod_{k=0}^{n-1}(1-aq^k)$. These polynomials are known to be orthogonal (but they are not orthonormal) on $[-1,1]$ with respect to a discrete measure $\sigma$ on the real line, concentrated at the two-sided sequence $\{\pm q^j\}_{n=0}^\infty$ with masses $\gamma_j \stackrel{\text{def}}{=} \sigma\{\pm q^j\} > 0$ for $j = 0, 1, \ldots$. The corresponding three-term recurrence is given by

$$H_{n+1}(x;q) = xH_n(x;q) - q^{n-1}(1-q^n)H_{n-1}(x;q), \quad H_0(x;q) = 1, \;\; H_1(x;q) = x, \quad n \in \mathbb{N},$$

(cf. [12, Exercise 7.38, p. 193]). Going over to the unit circle and denoting the corresponding measure by $\nu$, we have by Geronimus' formulas

$$\Phi_{2n-1}(\nu, 0) = 0 \quad \text{and} \quad \Phi_{2n}(\nu, 0) = 2q^n - 1, \qquad n \in \mathbb{N},$$

(cf. [16, Theorem 31.1, p. 67]). The measure $\nu$, which can be viewed as a measure on $[-\pi, \pi)$, is symmetric and it is concentrated at the set $\{\pm\theta_j^\pm\}_{j=0}^\infty$ where $\cos\theta_j^\pm = \pm q^j$ and $\nu\{\pm\theta_j^\pm\} = \gamma_j$ for $j = 0, 1, \ldots$. We have $\theta_j^+ + \theta_j^- = \pi$. The condition $\Phi_{2n-1}(\nu, 0) = 0$, $n \in \mathbb{N}$, means that the measure $\nu$ is "sieved". For the measure $\mu$ concentrated at the points $\{\tau_j^\pm\}$ where $\tau_j^+ \stackrel{\text{def}}{=} 2\theta_j^+$ and $\tau_j^- \stackrel{\text{def}}{=} 2\theta_j^- - 2\pi$ with masses $\mu\{\tau_j^\pm\} \stackrel{\text{def}}{=} \gamma_j$, we have $d\nu(\theta) = d\mu(2\theta)$. The support of the measure $\mu$ has only one limit point at $-\pi$, that is, the derived set of $\text{supp}(\mu)$ with $\mu$ viewed as a measure on $\mathbb{T}$, consists of $\{-1\}$, and corresponding recurrence coefficients $\Phi_n(\mu, 0)$ are equal to $2q^n - 1$ for $n \in \mathbb{N}$ (cf. [46, formula (7), p. 2] applied with $d = 1/2$).

## 4. Perturbation analysis

Given $a \in \mathbb{C}$ with $0 < |a| < 1$, consider

$$\widetilde{\varphi}_n(\mu, z) \stackrel{\text{def}}{=} \frac{\varphi_n(\mu, z)}{\kappa_n(\mu)(1-|a|^2)^{n/2}} \qquad \text{and} \qquad \widetilde{\psi}_n(\mu, z) \stackrel{\text{def}}{=} \frac{\psi_n(\mu, z)}{\kappa_n(\mu)(1-|a|^2)^{n/2}}. \qquad (32)$$

Then $\widetilde{\varphi}_n$ and $\widetilde{\psi}_n$ have the same leading coefficients as the comparison polynomials $\widehat{\varphi}_n$ and $\widehat{\psi}_n$ which are associated with the reflection coefficients $\{a\}_{n\in\mathbb{N}}$ and $\{-a\}_{n\in\mathbb{N}}$, respectively (cf. Section 2).

**Lemma 11.** *Given a probability measure $\mu$ on $\mathbb{T}$ having an infinite support, and $a \in \mathbb{C}$ with $0 < |a| < 1$, let $E_n(\mu, a)$ be defined by*

$$E_n(\mu, a, z) \stackrel{\text{def}}{=} \frac{1}{\sqrt{1-|a|^2}} \begin{pmatrix} 0 & \overline{\Phi_n(\mu,0)} - \overline{a} \\ z\,\overline{\Phi_n(\mu,0)} - \overline{a} & 0 \end{pmatrix}, \qquad n \in \mathbb{N}. \qquad (33)$$

*Then, for $n \in \mathbb{N}$,*

$$(1-|a|^2)^{n/2} \begin{pmatrix} z & a \\ z\overline{a} & 1 \end{pmatrix}^{-n} \begin{pmatrix} \widetilde{\varphi}_n(\mu,z) & \widetilde{\psi}_n(\mu,z) \\ \widetilde{\varphi}_n^*(\mu,z) & -\widetilde{\psi}_n^*(\mu,z) \end{pmatrix} = \quad (34)$$

$$\begin{pmatrix} 1 & 1 \\ 1 & -1 \end{pmatrix} + \sum_{k=0}^{n-1} (1-|a|^2)^{(k+1)/2} \begin{pmatrix} z & a \\ z\overline{a} & 1 \end{pmatrix}^{-k-1} E_{k+1}(\mu,a,z) \begin{pmatrix} \widetilde{\varphi}_k(\mu,z) & \widetilde{\psi}_k(\mu,z) \\ \widetilde{\varphi}_k^*(\mu,z) & -\widetilde{\psi}_k^*(\mu,z) \end{pmatrix},$$

*or, equivalently,*

$$\begin{pmatrix} \widetilde{\varphi}_n(\mu,z) & \widetilde{\psi}_n(\mu,z) \\ \widetilde{\varphi}_n^*(\mu,z) & -\widetilde{\psi}_n^*(\mu,z) \end{pmatrix} = \begin{pmatrix} \widehat{\varphi}_n(z) & \widehat{\psi}_n(z) \\ \widehat{\varphi}_n^*(z) & -\widehat{\psi}_n^*(z) \end{pmatrix} + \quad (35)$$

$$\frac{1}{2}\sum_{k=0}^{n-1} \begin{pmatrix} \widehat{\varphi}_{n-k-1}(z) & \widehat{\psi}_{n-k-1}(z) \\ \widehat{\varphi}_{n-k-1}^*(z) & -\widehat{\psi}_{n-k-1}^*(z) \end{pmatrix} \begin{pmatrix} 1 & 1 \\ 1 & -1 \end{pmatrix} E_{k+1}(\mu,a,z) \begin{pmatrix} \widetilde{\varphi}_k(\mu,z) & \widetilde{\psi}_k(\mu,z) \\ \widetilde{\varphi}_k^*(\mu,z) & -\widetilde{\psi}_k^*(\mu,z) \end{pmatrix}$$

*where $\{\widehat{\varphi}_n\}_{n\in\mathbb{N}}$ and $\{\widehat{\psi}_n\}_{n\in\mathbb{N}}$ are the orthogonal polynomials corresponding to the reflection coefficients $\{a\}_{n\in\mathbb{N}}$ and $\{-a\}_{n\in\mathbb{N}}$, respectively, that is,*

$$\begin{pmatrix} \widehat{\varphi}_n(z) & \widehat{\psi}_n(z) \\ \widehat{\varphi}_n^*(z) & -\widehat{\psi}_n^*(z) \end{pmatrix} = (1-|a|^2)^{-n/2} \begin{pmatrix} z & a \\ z\overline{a} & 1 \end{pmatrix}^n \begin{pmatrix} 1 & 1 \\ 1 & -1 \end{pmatrix} \quad (36)$$

*(cf. (3) and (15)).*

*Proof of Lemma 11.* We rewrite the matrix recurrence (8) as

$$\begin{pmatrix} \widetilde{\varphi}_{k+1}(\mu,z) & \widetilde{\psi}_{k+1}(\mu,z) \\ \widetilde{\varphi}_{k+1}^*(\mu,z) & -\widetilde{\psi}_{k+1}^*(\mu,z) \end{pmatrix} =$$

$$\frac{1}{\sqrt{1-|a|^2}} \begin{pmatrix} z & a \\ z\overline{a} & 1 \end{pmatrix} \begin{pmatrix} \widetilde{\varphi}_k(\mu,z) & \widetilde{\psi}_k(\mu,z) \\ \widetilde{\varphi}_k^*(\mu,z) & -\widetilde{\psi}_k^*(\mu,z) \end{pmatrix} + E_{k+1}(\mu,a,z) \begin{pmatrix} \widetilde{\varphi}_k(\mu,z) & \widetilde{\psi}_k(\mu,z) \\ \widetilde{\varphi}_k^*(\mu,z) & -\widetilde{\psi}_k^*(\mu,z) \end{pmatrix}.$$

(cf. (33)). Multiplying this from the left by $(1-|a|^2)^{(k+1)/2} \begin{pmatrix} z & a \\ z\overline{a} & 1 \end{pmatrix}^{-k-1}$ and summing over $k=0,1,\ldots,n-1$, yields (34). Formula (35) follows directly from (34) and (36). □

**Theorem 12.** *Let $\mu$ be a probability measure on $\mathbb{T}$ having an infinite support. Given $a \in \mathbb{C}$ with $0 < |a| < 1$, let $\alpha \in (0,\pi)$ be defined by $\cos\frac{\alpha}{2} \overset{\text{def}}{=} \sqrt{1-|a|^2}$. If the reflection coefficients $\{\Phi_k(\mu,0)\}_{k\in\mathbb{N}}$ of the corresponding orthogonal polynomials satisfy*

$$\sum_{k=1}^{\infty} |\tau^k \Phi_k(\mu,0) - a| < \infty, \qquad a \in \mathbb{C},\ 0 < |a| < 1,\ \tau \in \mathbb{T}, \quad (37)$$

*then the measure $\mu$ is absolutely continuous on the open circular arc[11] $\tau\Delta_\alpha^o$, and for every closed circular subarc $\mathcal{E} \subset \tau\Delta_\alpha^o$ we have $1/\mu' \in L^\infty(\mathcal{E})$. If*

$$\sum_{k=1}^{\infty} \log k\, |\tau^k \Phi_k(\mu,0) - a| < \infty, \qquad a \in \mathbb{C},\ 0 < |a| < 1,\ \tau \in \mathbb{T}, \quad (38)$$

---

[11] Recall that the circular arcs $\Delta_\alpha$ (closed) and $\Delta_\alpha^o$ (open) have been defined in (13) and (14). The rotation of a circular arc $\Delta$ by $\arg \tau$ for $\tau \in \mathbb{T}$ is denoted by $\tau\Delta$.

then $\mu$ satisfies Szegő's condition on $\tau \Delta_\alpha$, that is, $v_{\alpha,\tau} \log \mu' \in L^1(\tau \Delta_\alpha)$ where $v_{\alpha,\tau}(\theta) \stackrel{\text{def}}{=}$ $\sin \frac{\theta - \arg \tau}{2} |\cos \alpha - \cos(\theta - \arg \tau)|^{-\frac{1}{2}}$ (cf. (22)). If there is a constant $C > 0$ such that

$$\sum_{k=1}^n k |\tau^k \Phi_k(\mu, 0) - a| < C \log(n+1), \qquad a \in \mathbb{C}, \; 0 < |a| < 1, \; \tau \in \mathbb{T}, \; n \in \mathbb{N}, \qquad (39)$$

then there exist two constants $D > 0$ and $\gamma > 0$ such that $\mu'(z) \geq D |\cos \alpha - \cos(\theta - \arg \tau)|^\gamma$ for almost every $z = e^{i\theta} \in \tau \Delta_\alpha$. Moreover, if

$$\sum_{k=1}^\infty k |\tau^k \Phi_k(\mu, 0) - a| < \infty, \qquad a \in \mathbb{C}, \; 0 < |a| < 1, \; \tau \in \mathbb{T}, \qquad (40)$$

then there exists a constant $D > 0$ such that $\mu'(z) \geq D |\cos \alpha - \cos(\theta - \arg \tau)|$ for almost every $z = e^{i\theta} \in \tau \Delta_\alpha$.

*Proof of Theorem 12.* First, observe that it is sufficient to prove the theorem for $\tau = 1$. To see this, note that a composition of the measure $\mu$ with the rotation $\tau^{-1}$ leads to the orthogonal polynomials $\{\varphi_n(\mu \circ \tau^{-1}, z) \equiv \tau^n \varphi_n(\mu, \tau^{-1} z)\}_{n=0}^\infty$ so that the zeros of the orthogonal polynomials associated with $\mu \circ \tau^{-1}$ are $\tau$ times the zeros of the orthogonal polynomials corresponding to $\mu$. Hence the reflection coefficients, which are the products of such zeros, satisfy $\Phi_n(\mu \circ \tau^{-1}, 0) = \tau^n \Phi_n(\mu, 0)$ for $n \in \mathbb{N}$. Therefore, we can assume without loss of generality that $\tau = 1$.

Second, observe that, by (3), either of (37)–(40) imply

$$0 < \lim_{n \to \infty} \kappa_n(\mu)(1 - |a|^2)^{n/2} < \infty. \qquad (41)$$

Hence, $\sup_{n \in \mathbb{N}} \max_{z \in \mathbb{T}} |\widetilde{\varphi}_n(\mu, z)/\varphi_n(\mu, z)|^{\pm 1}$ and $\sup_{n \in \mathbb{N}} \max_{z \in \mathbb{T}} |\widetilde{\psi}_n(\mu, z)/\psi_n(\mu, z)|^{\pm 1}$ are both finite (cf. (32)).

We will use Lemmas 1 and 2. In what follows, we fix a norm for $2 \times 2$ matrices and define $a$ by the appropriate condition (37)–(40).

**Case when (37) holds.** Pick a closed circular subarc $\mathcal{E} \subset \Delta_\alpha^o$. Then, by (25) and (35), there are two positive constants $C_1$ and $C_2$ independent of $n$ such that

$$\max_{z \in \mathcal{E}} \left\| \begin{pmatrix} \widetilde{\varphi}_n(\mu, z) & \widetilde{\psi}_n(\mu, z) \\ \widetilde{\varphi}_n^*(\mu, z) & -\widetilde{\psi}_n^*(\mu, z) \end{pmatrix} \right\| \leq$$
$$C_1 + C_2 \sum_{k=0}^{n-1} \max_{z \in \mathcal{E}} \|E_{k+1}(\mu, a, z)\| \max_{z \in \mathcal{E}} \left\| \begin{pmatrix} \widetilde{\varphi}_k(\mu, z) & \widetilde{\psi}_k(\mu, z) \\ \widetilde{\varphi}_k^*(\mu, z) & -\widetilde{\psi}_k^*(\mu, z) \end{pmatrix} \right\|, \qquad n \in \mathbb{N},$$

so that by Gronwall's inequality (cf. [45, Lemma, p. 440] and the references therein),

$$\max_{z \in \mathcal{E}} \left\| \begin{pmatrix} \widetilde{\varphi}_n(\mu, z) & \widetilde{\psi}_n(\mu, z) \\ \widetilde{\varphi}_n^*(\mu, z) & -\widetilde{\psi}_n^*(\mu, z) \end{pmatrix} \right\| \leq C_1 \exp\left( C_2 \sum_{k=0}^{n-1} \max_{z \in \mathcal{E}} \|E_{k+1}(\mu, a, z)\| \right), \qquad n \in \mathbb{N}.$$

Thus, since $\max_{z\in\mathcal{E}} \|E_{k+1}(\mu,a,z)\| \leq C_3|\Phi_{k+1}(\mu,0) - a|$ (cf. (33)) with an appropriate positive constant $C_3$ independent of $k$,

$$\sup_{n\in\mathbb{N}} \max_{z\in\mathcal{E}} \left\| \begin{pmatrix} \widetilde{\varphi}_n(\mu,z) & \widetilde{\psi}_n(\mu,z) \\ \widetilde{\varphi}_n^*(\mu,z) & -\widetilde{\psi}_n^*(\mu,z) \end{pmatrix} \right\| < \infty, \qquad (42)$$

and, therefore, by Lemmas 1 and 2, $\mu \in AC(\mathcal{E})$ and $1/\mu' \in L^\infty(\mathcal{E})$, respectively (cf. (32) and (41)).

**Case when (38) holds.** The basic idea goes back to [29, Theorem 6, p. 381] and [7, Theorem 3, p. 355]; it is based on estimates of $\log^+ 1/\mu'$ in terms of $\limsup_{n\to\infty} \log^+ |\varphi_n(\mu)|$ (cf. (12)).

Let

$$S_n(\mu,a,z) \stackrel{\text{def}}{=} 1 + \sum_{k=0}^{n-1} \|E_{k+1}(\mu,a,z)\| \left\| \begin{pmatrix} \widetilde{\varphi}_k(\mu,z) & \widetilde{\psi}_k(\mu,z) \\ \widetilde{\varphi}_k^*(\mu,z) & -\widetilde{\psi}_k^*(\mu,z) \end{pmatrix} \right\|, \qquad n \in \mathbb{Z}^+, \quad (43)$$

(cf. (33)). Then, by (24) and (35),

$$S_n(\mu,a,z) = S_{n-1}(\mu,a,z) + \|E_n(\mu,a,z)\| \left\| \begin{pmatrix} \widetilde{\varphi}_{n-1}(\mu,z) & \widetilde{\psi}_{n-1}(\mu,z) \\ \widetilde{\varphi}_{n-1}^*(\mu,z) & -\widetilde{\psi}_{n-1}^*(\mu,z) \end{pmatrix} \right\|$$

$$\leq S_{n-1}(\mu,a,z) \left[1 + C_1^* \min(n, v_\alpha(\theta)) \|E_n(\mu,a,z)\|\right], \qquad z = e^{i\theta},\ z \in \Delta_\alpha,\ n \in \mathbb{N},$$

with an appropriate constant $C_1^* > 0$, so that, by induction,

$$0 \leq \log S_n(\mu,a,z) \leq \log S_1(\mu,a,z) + C_1^* \sum_{k=2}^{n} \min(k, v_\alpha(\theta)) \|E_k(\mu,a,z)\|,$$

$$z = e^{i\theta},\ z \in \Delta_\alpha,\ n \in \mathbb{N}. \qquad (44)$$

Multiply (44) by $v_\alpha$ and integrate the resulting inequality over $\Delta_\alpha$. Taking (38) and (33) into consideration, and using $\sup_{n\in\mathbb{N}} \left[\frac{1}{\log(n+1)} \int_{\Delta_\alpha} v_\alpha \min(n, v_\alpha)\right] < \infty$, simple computation yields $\sup_{n\in\mathbb{N}} \int_{\Delta_\alpha} v_\alpha \log S_n(\mu,a) < \infty$, and then, by Lebesgue's *monotone convergence theorem*, $\sup_{n\in\mathbb{N}} v_\alpha \log S_n(\mu,a) \in L^1(\Delta_\alpha)$. Hence, by (24), (35), and (43), $\sup_{n\in\mathbb{N}} v_\alpha \log^+ |\varphi_n(\mu)| \in L^1(\Delta_\alpha)$. Now use Lemma 2 to conclude that $v_\alpha \log^+(1/\mu') \in L^1(\Delta_\alpha)$. The integrability of $v_\alpha \log^-(1/\mu')$ follows from Jensen's (AGM) inequality.

**Case when (39) holds.** We will adapt the arguments used in [28]. Our goal is to show that there exists a constant $\gamma > 0$ such that

$$D^{-\frac{1}{2}} \stackrel{\text{def}}{=} \sqrt{2\pi} \sup_{n\in\mathbb{N}} \max_{z\in\Delta_\alpha} |\cos\alpha - \cos\theta|^{\frac{\gamma}{2}} |\varphi_n(\mu,z)| < \infty, \qquad z = e^{i\theta}, \qquad (45)$$

and then, by Lemma 2, $\mu'(z) > D |\cos\alpha - \cos\theta|^\gamma$ for almost every $z = e^{i\theta} \in \Delta_\alpha$ what needs to be proved.

By (24) and (35), there are two positive constants $C_4$ and $C_5$ independent of $n$ such that

$$\frac{1}{n} \max_{z \in \Delta_\alpha} \left\| \begin{pmatrix} \widetilde{\varphi}_n(\mu, z) & \widetilde{\psi}_n(\mu, z) \\ \widetilde{\varphi}_n^*(\mu, z) & -\widetilde{\psi}_n^*(\mu, z) \end{pmatrix} \right\| \leq$$
$$C_4 + C_5 \sum_{k=0}^{n-1} k \max_{z \in \Delta_\alpha} \|E_{k+1}(\mu, a, z)\| \frac{1}{k} \max_{z \in \Delta_\alpha} \left\| \begin{pmatrix} \widetilde{\varphi}_k(\mu, z) & \widetilde{\psi}_k(\mu, z) \\ \widetilde{\varphi}_k^*(\mu, z) & -\widetilde{\psi}_k^*(\mu, z) \end{pmatrix} \right\|, \quad n \in \mathbb{N},$$

so that by Gronwall's inequality (cf. [45, Lemma, p. 440]),

$$\frac{1}{n} \max_{z \in \Delta_\alpha} \left\| \begin{pmatrix} \widetilde{\varphi}_n(\mu, z) & \widetilde{\psi}_n(\mu, z) \\ \widetilde{\varphi}_n^*(\mu, z) & -\widetilde{\psi}_n^*(\mu, z) \end{pmatrix} \right\| \leq C_4 \exp\left( C_5 \sum_{k=0}^{n-1} k \max_{z \in \Delta_\alpha} \|E_{k+1}(\mu, a, z)\| \right), \quad n \in \mathbb{N}. \tag{46}$$

Thus, since $\|E_{k+1}(\mu, a, z)\| \leq C_3 |\Phi_{k+1}(0) - a|$ (cf. (33)) with an appropriate positive constant $C_3$ independent of $k$, assumption (39) guarantees the existence of a constant $\gamma^* > 0$ such that

$$\sup_{n \in \mathbb{N}} \max_{z \in \Delta_\alpha} n^{-\gamma^*} \left\| \begin{pmatrix} \widetilde{\varphi}_n(\mu, z) & \widetilde{\psi}_n(\mu, z) \\ \widetilde{\varphi}_n^*(\mu, z) & -\widetilde{\psi}_n^*(\mu, z) \end{pmatrix} \right\| < \infty. \tag{47}$$

In the following argument, we will use the inequality

$$|\cos \alpha - \cos \theta|^{\frac{1}{2}} \left\| \begin{pmatrix} \widetilde{\varphi}_n(\mu, z) & \widetilde{\psi}_n(\mu, z) \\ \widetilde{\varphi}_n^*(\mu, z) & -\widetilde{\psi}_n^*(\mu, z) \end{pmatrix} \right\| \leq \tag{48}$$
$$C_6 + C_6 \sum_{k=0}^{n-1} |\Phi_{k+1}(0) - a| \left\| \begin{pmatrix} \widetilde{\varphi}_k(\mu, z) & \widetilde{\psi}_k(\mu, z) \\ \widetilde{\varphi}_k^*(\mu, z) & -\widetilde{\psi}_k^*(\mu, z) \end{pmatrix} \right\|,$$
$$z = e^{i\theta}, \ z \in \Delta_\alpha, \ n \in \mathbb{N},$$

where $C_6$ is an appropriate positive constant. This is an immediate consequence of (24) and (35) (cf. (33)). Also note that by (39)

$$\sum_{k=1}^{\infty} k^\sigma |\Phi_k(\mu, 0) - a| < \infty, \quad 0 \leq \sigma < 1. \tag{49}$$

If $\gamma^* < 1$ in (47) then, by (48) and (49), inequality (45) holds with $\gamma = 1$ and we are done. Otherwise, using (47) on the right–hand side of (48) and applying (49) with $\sigma = 1/2$ yields

$$\sup_{n \in \mathbb{N}} \max_{z \in \Delta_\alpha} n^{-\gamma^* + \frac{1}{2}} |\cos \alpha - \cos \theta|^{\frac{1}{2}} \left\| \begin{pmatrix} \widetilde{\varphi}_n(\mu, z) & \widetilde{\psi}_n(\mu, z) \\ \widetilde{\varphi}_n^*(\mu, z) & -\widetilde{\psi}_n^*(\mu, z) \end{pmatrix} \right\| < \infty. \tag{50}$$

If $\gamma^* < 3/2$ in (50) then again, by (48) and (49), inequality (45) holds with $\gamma = 2$ and we are done. If $\gamma^* \geq 3/2$ then we continue this procedure each step of which lowers $\gamma^*$ in the analogues of (47) and (50) by $1/2$ and, simultaneously, raises the value of $\gamma$ in (45) by 1.

After no more than $[2\gamma^*] - 1$ steps[12] this procedure terminates and then (45) holds with $\gamma \leq [2\gamma^*]$.

**Case when** (40) **holds.** The inequality $\mu'(z) \geq D \, |\cos\alpha - \cos\theta|$ for almost every $z = e^{i\theta} \in \Delta_\alpha$ can be proved along the lines of the case when (39) holds. Namely, by (46) and (48) (cf. (33) and (40)),

$$D^{-\frac{1}{2}} \stackrel{\text{def}}{=} \sqrt{2\pi} \sup_{n \in \mathbb{N}} \max_{z \in \Delta_\alpha} |\cos\alpha - \cos\theta|^{\frac{1}{2}} |\varphi_n(\mu, z)| < \infty, \qquad z = e^{i\theta},$$

which is the analogue of (45), and then we just repeat the argument used to obtain a lower estimate for $\mu'$ from (45), that is, use Lemma 2.
□

*Remark 13.* F. Peherstorfer and R. Steinbauer used a different approach in [35] to study the asymptotic behavior of orthogonal polynomials with asymptotically periodic recurrence coefficients.

The following theorem is a useful summary of two inequalities proved in the proof of Theorem 12.

**Theorem 14.** *Let $\mu$ be a probability measure on $\mathbb{T}$ having an infinite support. Given $a \in \mathbb{C}$ with $0 < |a| < 1$, let $\alpha \in (0, \pi)$ be defined by $\cos\frac{\alpha}{2} \stackrel{\text{def}}{=} \sqrt{1 - |a|^2}$. If the reflection coefficients $\{\Phi_k(\mu, 0)\}_{k \in \mathbb{N}}$ of the corresponding orthogonal polynomials satisfy*

$$\sum_{k=1}^{\infty} |\tau^k \Phi_k(\mu, 0) - a| < \infty, \qquad a \in \mathbb{C}, \ 0 < |a| < 1, \ \tau \in \mathbb{T}, \tag{51}$$

*then for every closed circular subarc[13] $\mathcal{E} \subset \tau \Delta_\alpha^o$ we have*

$$\sup_{n \in \mathbb{N}} \max_{z \in \mathcal{E}} |\varphi_n(\mu, z)| < \infty. \tag{52}$$

*If*

$$\sum_{k=1}^{\infty} k |\tau^k \Phi_k(\mu, 0) - a| < \infty, \qquad a \in \mathbb{C}, \ 0 < |a| < 1, \ \tau \in \mathbb{T}, \tag{53}$$

*then*

$$\sup_{n \in \mathbb{N}} \frac{\max_{z \in \tau \Delta_\alpha} |\varphi_n(\mu, z)|}{n} < \infty. \tag{54}$$

*Proof of Theorem 14.* Just like in the proof of Theorem 12, it is sufficient to prove the theorem for $\tau = 1$. In addition, again as in the proof of Theorem 12, in view of (32) and (41), either of the assumptions (51) and (53) imply

$$\sup_{n \in \mathbb{N}} \max_{z \in \mathbb{T}} \frac{|\varphi_n(\mu, z)|}{|\widetilde{\varphi}_n(\mu, z)|} < \infty$$

---

[12]Here $[\cdot]$ denotes the *integer part*.
[13]Recall that the circular arcs $\Delta_\alpha$ (closed) and $\Delta_\alpha^o$ (open) have been defined in (13) and (14). The rotation of a circular arc $\Delta$ by $\arg\tau$ for $\tau \in \mathbb{T}$ is denoted by $\tau\Delta$.

so that we need to estimate $|\widetilde{\varphi}_n(\mu, z)|$ only. If (51) holds then (52) follows from (42), whereas if (53) holds then (54) follows from (46) (cf. (33)). □

*Remark 15.* F. Peherstorfer and R. Steinbauer used a different approach in [36, Proposition 2.1] to prove (52) in Theorem 14.

## 5. Perturbation analysis (continued)

The condition $\sum_{k \in \mathbb{N}} |\Phi_k(\mu, 0) - a| < \infty$ can be weakened. To do this one needs to write the perturbation series in a multiplicative way, much as it is done in [34] and [31]. Clearly, by (8) and (32),

$$\begin{pmatrix} \widetilde{\varphi}_n(\mu, z) & \widetilde{\psi}_n(\mu, z) \\ \widetilde{\varphi}_n^*(\mu, z) & -\widetilde{\psi}_n^*(\mu, z) \end{pmatrix} = (\Omega(a, z) + E_n(\mu, a, z)) \begin{pmatrix} \widetilde{\varphi}_{n-1}(\mu, z) & \widetilde{\psi}_{n-1}(\mu, z) \\ \widetilde{\varphi}_{n-1}^*(\mu, z) & -\widetilde{\psi}_{n-1}^*(\mu, z) \end{pmatrix}, \quad n \in \mathbb{N},$$

where

$$\Omega(a, z) \stackrel{\text{def}}{=} \frac{1}{\sqrt{1 - |a|^2}} \begin{pmatrix} z & a \\ z\overline{a} & 1 \end{pmatrix}.$$

Writing

$$B_n(\mu, a, z) \stackrel{\text{def}}{=} (1 - |a|^2)^{n/2} \begin{pmatrix} z & a \\ z\overline{a} & 1 \end{pmatrix}^{-n} \begin{pmatrix} \widetilde{\varphi}_n(\mu, z) & \widetilde{\psi}_n(\mu, z) \\ \widetilde{\varphi}_n^*(\mu, z) & -\widetilde{\psi}_n^*(\mu, z) \end{pmatrix},$$

one obtains

$$B_n(\mu, a, z) = B_{n-1}(\mu, a, z) + \Omega^{-n}(a, z) E_n(\mu, a, z) \Omega^{n-1}(a, z) B_{n-1}(\mu, a, z)$$
$$= \left( I + \Omega^{-n}(a, z) E_n(\mu, a, z) \Omega^{n-1}(a, z) \right) B_{n-1}(\mu, a, z),$$

and then iterating this yields

$$B_n(\mu, a, z) = \prod_{k=1}^n \left( I + \Omega^{-k}(a, z) E_k(\mu, a, z) \Omega^{k-1}(a, z) \right) \begin{pmatrix} 1 & 1 \\ 1 & -1 \end{pmatrix}, \quad n \in \mathbb{N},$$

where the matrix product is meant to be taken from the right ($k = 1$) to the left ($k = n$). Observe that

$$\Omega^n(a, z) = \frac{1}{2} \begin{pmatrix} \widehat{\varphi}_n(z) & \widehat{\psi}_n(z) \\ \widehat{\varphi}_n^*(z) & -\widehat{\psi}_n^*(z) \end{pmatrix} \begin{pmatrix} 1 & 1 \\ 1 & -1 \end{pmatrix}$$

and

$$\Omega^{-n}(a, z) = \frac{1}{2z^n} \begin{pmatrix} 1 & 1 \\ 1 & -1 \end{pmatrix} \begin{pmatrix} \widehat{\psi}_n^*(z) & \widehat{\psi}_n(z) \\ \widehat{\varphi}_n^*(z) & -\widehat{\varphi}_n(z) \end{pmatrix}$$

(cf. (36)). We can thus formulate the following result

**Theorem 16.** *Suppose the infinite matrix product*

$$\prod_{k=1}^{\infty} \left(I + \Omega^{-k}(a,z) E_k(\mu, a, z) \Omega^{k-1}(a,z)\right) \tag{55}$$

*converges (conditionally) in the open circular arc $\Delta_\alpha^o$. Then*

$$\lim_{n\to\infty} (1-|a|^2)^{n/2} \begin{pmatrix} z & a \\ z\bar{a} & 1 \end{pmatrix}^{-n} \begin{pmatrix} \widetilde{\varphi}_n(\mu,z) & \widetilde{\psi}_n(\mu,z) \\ \widetilde{\varphi}_n^*(\mu,z) & -\widetilde{\psi}_n^*(\mu,z) \end{pmatrix}$$

*exists in $\Delta_\alpha^o$.*

The interesting and challenging problem is to describe the conditional convergence of (55) in terms of $\Phi_k(\mu,0) - a$. For instance, if $\sum_{k=1}^\infty [\Phi_k(\mu,0) - a]$ converges (conditionally) and $\sum_{k=1}^\infty |\Phi_k(\mu,0) - a|^2 < \infty$, then what additional conditions (if any) are needed to assure the convergence of (55)?

## 6. An example

**Example 17.** Given $\alpha \in [0,\pi)$, $\gamma > -1$, and $\delta > -1$, consider the absolutely continuous measure $\mu$ in $[0, 2\pi)$ with

$$\mu'(\theta) \stackrel{\text{def}}{=} \begin{cases} C\, |\cos\theta - \cos\alpha|^\gamma |\cos\tfrac{\theta}{2}|^\delta |\sin\tfrac{\theta}{2}| & \text{if } \theta \in (\alpha, 2\pi - \alpha), \\ 0 & \text{if } \theta \in [0, 2\pi) \setminus (\alpha, 2\pi - \alpha), \end{cases}$$

where the normalizing factor $C > 0$ is chosen in such a way that $\mu(\mathbb{T}) = 1$. Then the reflection coefficients for this measure satisfy

$$\Phi_n(\mu, 0) = \sin\tfrac{\alpha}{2} + (-1)^n \frac{\delta \cos^2\tfrac{\alpha}{2}}{2n} + \frac{\cos\tfrac{\alpha}{2} \cot\tfrac{\alpha}{2}}{16n^2}(-2 - \delta^2 + 8\gamma^2 + \delta^2 \cos\alpha)$$
$$- (-1)^n \frac{\delta \cos^2\tfrac{\alpha}{2}}{4n^2} (1 + \delta + 2\gamma - \sin\tfrac{\alpha}{2}) + O(1/n^3), \qquad n \in \mathbb{N}. \tag{56}$$

In particular, for this measure $\mu$ and for $a \stackrel{\text{def}}{=} \sin\tfrac{\alpha}{2}$, formulas (37) and (38) hold if and only if $\delta = 0$, whereas formulas (39) and (40) hold if and only if $\delta = 0$ and $\gamma^2 = 1/4$. This is in complete agreement with Theorem 12. The case $\delta = 0$ and $\gamma^2 = 1/4$ corresponds to the Lee–Yang weight function associated with the one dimensional Ising model (cf. [24, formula (45), p. 416]).

*Proof of* (56). First we make a transition from the unit circle to the real line (cf. [42, § 11.5, p. 294]). The relevant weight function on the real line is $(\cos\alpha - x)^\gamma (1+x)^{\frac{\delta-1}{2}}$ supported in $[-1, \cos\alpha]$. If we map the interval $[-1, \cos\alpha]$ to $[-1, 1]$ by the affine transformation $x \mapsto \frac{2x + 1 - \cos\alpha}{1 + \cos\alpha}$, then the orthogonal polynomials will be the Jacobi polynomials $P_n^{(\gamma, \frac{\delta-1}{2})}$, and, hence, the orthogonal polynomials on $[-1, \cos \alpha]$ are $P_n^{(\gamma, \frac{\delta-1}{2})}\left(\frac{2x+1-\cos \alpha}{1+\cos \alpha}\right)$. The monic orthogonal polynomials are given by

$$P_n(x) \stackrel{\text{def}}{=} (1+\cos\alpha)^n \binom{2n+\gamma+\frac{\delta-1}{2}}{n}^{-1} P_n^{(\gamma, \frac{\delta-1}{2})}\left(\frac{2x+1-\cos\alpha}{1+\cos\alpha}\right), \qquad n \in \mathbb{N}.$$

The reflection coefficients $\Phi_n(\mu, 0)$ can be obtained from these monic polynomials by

$$\Phi_{2n}(\mu, 0) \equiv R_n(1) - R_n(-1) - 1 \quad \text{and} \quad \Phi_{2n+1}(\mu, 0) \equiv \frac{R_n(1) + R_n(-1)}{R_n(1) - R_n(-1)}$$

where $R_n \stackrel{\text{def}}{=} P_{n+1}/P_n$ (cf. [16, Theorem 31.1, p. 67] or [15, Theorem X, p. 758]). In our case, we have

$$R_n(1) = (1+\cos\alpha)\frac{(n+1)\left(n+\gamma+\frac{\delta}{2}+\frac{1}{2}\right)}{\left(2n+\gamma+\frac{\delta}{2}+\frac{1}{2}\right)\left(2n+\gamma+\frac{\delta}{2}+\frac{3}{2}\right)} \frac{P_{n+1}^{(\gamma, \frac{\delta-1}{2})}\left(\frac{3-\cos\alpha}{1+\cos\alpha}\right)}{P_n^{(\gamma, \frac{\delta-1}{2})}\left(\frac{3-\cos\alpha}{1+\cos\alpha}\right)}$$

and

$$R_n(-1) = -(1+\cos\alpha)\frac{\left(n+\frac{\delta}{2}+\frac{1}{2}\right)\left(n+\gamma+\frac{\delta}{2}+\frac{1}{2}\right)}{\left(2n+\gamma+\frac{\delta}{2}+\frac{1}{2}\right)\left(2n+\gamma+\frac{\delta}{2}+\frac{3}{2}\right)}$$

where we used the fact that $P_n^{(\gamma, \frac{\delta-1}{2})}(-1) = (-1)^n \binom{n+\frac{\delta-1}{2}}{n}$ (cf. [42, formula (4.1.4), p. 59]). We have

$$\frac{(n+1)(n+\gamma+\frac{\delta}{2}+\frac{1}{2})}{(2n+\gamma+\frac{\delta}{2}+\frac{1}{2})(2n+\gamma+\frac{\delta}{2}+\frac{3}{2})} = \frac{1}{4} + \frac{\delta}{8n} + \frac{1 - 4\delta - \delta^2 - 4\delta\gamma - 4\gamma^2}{64n^2} + O\left(n^{-3}\right),$$

($n \in \mathbb{N}$) which leads to

$$R_n(-1) = -\frac{1+\cos\alpha}{4}\left(1 + \frac{\delta}{2n} + \frac{1 - 4\delta - \delta^2 - 4\delta\gamma - 4\gamma^2}{16n^2} + O\left(n^{-3}\right)\right), \qquad n \in \mathbb{N}.$$

For estimating $R_n(1)$, we need an asymptotic expression for the ratio of two consecutive Jacobi polynomials at the point $(3-\cos\alpha)/(1+\cos\alpha) > 1$.

**Lemma 18.** *Given $a > -1$, $b > -1$, and $x > 1$, the Jacobi polynomials $\{P_n^{(a,b)}\}_{n \in \mathbb{N}}$ normalized by $P_n^{(a,b)}(1) = \binom{n+a}{n}$ satisfy the ratio asymptotics*

$$\frac{(n+1)(n+a+b+1)}{(2n+a+b+2)(2n+a+b+1)} \frac{P_{n+1}^{(a,b)}(x)}{P_n^{(a,b)}(x)}$$

$$= \frac{x+\sqrt{x^2-1}}{4}\left(1 - \frac{A(x,a,b)}{2n^2} + O\left(n^{-3}\right)\right), \qquad n \in \mathbb{N},$$

*where*

$$A(x, a, b) \stackrel{\text{def}}{=} \frac{b^2}{x+1+\sqrt{x^2-1}} - \frac{a^2}{x-1+\sqrt{x^2-1}} + \frac{1}{4}\frac{1}{\sqrt{x^2-1}\,(x+\sqrt{x^2-1})}\,.$$

*Proof of Lemma 18.* By the asymptotic expansion [9, formula (2.13)]

$$P_n^{(a,b)}(x) = \frac{\Gamma(2n+a+b+1)}{2^{2n+(a+b+1)/2}\Gamma(n+1)\Gamma(n+a+b+1)}$$
$$\times \frac{(x+\sqrt{x^2-1})^{n+(a+b+1)/2}}{(x-1)^{(2a+1)/4}(x+1)^{(2b+1)/4}} \sum_{k=0}^{\infty} \frac{a_k(y)}{(2n+a+b+1)^k}\,, \qquad n \in \mathbb{N},$$

where $x = \cosh 2y \notin [-1, 1]$, a straightforward (but rather tedious) computation readily leads to the desired result taking into account that $a_0 \equiv 1$ and $A(x, a, b) \equiv a_1(y)$. The function $a_1$ can be computed explicitly, by

$$a_1(y) \equiv -\frac{1}{2}\int_y^\infty \left[\frac{2a^2}{\cosh(2t)-1} - \frac{2b^2}{\cosh(2t)+1} - \frac{1}{\sinh^2(2t)}\right] dt$$

(see [9, formulas (2.8) and (2.12)]). Simple calculus gives

$$\int_y^\infty \frac{dt}{\cosh(2t)\pm 1} = \frac{1}{e^{2y}\pm 1} = \frac{1}{x\pm 1 + \sqrt{x^2-1}}$$

and

$$\int_y^\infty \frac{dt}{\sinh^2(2t)} = \frac{1}{e^{4y}-1} = \frac{1}{2\sqrt{x^2-1}\,(x+\sqrt{x^2-1})}$$

which yields the desired expression for $A(x, a, b)$. □

Continuing with our computation, simple trigonometry shows[14]

$$x = (3-\cos\alpha)/(1+\cos\alpha) \quad \Longrightarrow \quad x+\sqrt{x^2-1} = \frac{2(1+\sin\frac{\alpha}{2})^2}{1+\cos\alpha}\,,$$

so that, by the lemma,

$$R_n(1) = \frac{(1+\sin\frac{\alpha}{2})^2}{2}\left(1 - \frac{A\left(x, \gamma, \frac{\delta-1}{2}\right)}{2n^2} + O\left(n^{-3}\right)\right).$$

---

[14]Curiously, Mathematica™ *failed* to prove this even though it was used successfully to perform and/or to double check much of the computation in this section.

Combining the formulas for $R_n(-1)$ and $R_n(1)$ gives

$$\Phi_{2n}(\mu,0) \equiv R_n(1) - R_n(-1) - 1 = \sin\tfrac{\alpha}{2} + \frac{\delta\cos^2\tfrac{\alpha}{2}}{4n}$$
$$+ \frac{\cos\tfrac{\alpha}{2}\cot\tfrac{\alpha}{2}}{64n^2}\left(-2 + 2\delta - \delta^2 + 8\gamma^2 - 2\delta\cos\alpha + \delta^2\cos\alpha\right.$$
$$\left. - 4\delta\sin\tfrac{\alpha}{2} - 4\delta^2\sin\tfrac{\alpha}{2} - 8\delta\gamma\sin\tfrac{\alpha}{2}\right) + O\left(n^{-3}\right), \qquad n \in \mathbb{N},$$

and

$$\Phi_{2n+1}(\mu,0) \equiv \frac{R_n(1) + R_n(-1)}{R_n(1) - R_n(-1)} = \sin\tfrac{\alpha}{2} - \frac{\delta\cos^2\tfrac{\alpha}{2}}{4n}$$
$$+ \frac{\cos\tfrac{\alpha}{2}\cot\tfrac{\alpha}{2}}{64n^2}\left(-2 - 2\delta - \delta^2 + 8\gamma^2 + 2\delta\cos\alpha + \delta^2\cos\alpha\right.$$
$$\left. + 12\delta\sin\tfrac{\alpha}{2} + 4\delta^2\sin\tfrac{\alpha}{2} + 8\delta\gamma\sin\tfrac{\alpha}{2}\right) + O\left(n^{-3}\right), \qquad n \in \mathbb{N},$$

that is,

$$\Phi_{2n+1}(\mu,0) = \sin\tfrac{\alpha}{2} - \frac{\delta\cos^2\tfrac{\alpha}{2}}{2(2n+1)}$$
$$+ \frac{\cos\tfrac{\alpha}{2}\cot\tfrac{\alpha}{2}}{16(2n+1)^2}\left(-2 - 2\delta - \delta^2 + 8\gamma^2 + 2\delta\cos\alpha + \delta^2\cos\alpha\right.$$
$$\left. + 4\delta\sin\tfrac{\alpha}{2} + 4\delta^2\sin\tfrac{\alpha}{2} + 8\delta\gamma\sin\tfrac{\alpha}{2}\right) + O\left(n^{-3}\right), \qquad n \in \mathbb{N},$$

and, therefore, we have proved (56). □

(L.G.) Mathematical Division, Institute for Low Temperature Physics and Engineering, 47 Lenin Avenue, Kharkov 310164, Ukraine
 *E-mail address*: golinskii@ilt.kharkov.ua

(P.N.) Department of Mathematics, Ohio State University, 231 West 18th Avenue, Columbus, Ohio 43210-1174, U.S.A.
 *E-mail address*: nevai@math.ohio-state.edu — *WWW home*: http://www.math.ohio-state.edu/~nevai

(W.V.A.) Department of Mathematics, Katholieke Universiteit Leuven, Celestijnenlaan 200 B, B–3001 Heverlee (Leuven), BELGIUM
 *E-mail address*: Walter.VanAssche@wis.KULeuven.ac.be